\documentclass[a4paper,twoside,12pt]{amsart}

\usepackage[utf8]{inputenc}
\usepackage[T1]{fontenc}
\usepackage[english]{babel}
\usepackage{tgpagella}
\usepackage{hyperref}
\usepackage{amscd} 
\usepackage{csquotes}
\usepackage{fancyhdr}
\usepackage{multirow}

\usepackage{amsthm, amsfonts, amssymb, amsmath, latexsym, enumerate,array}
\usepackage{amscd} 
\usepackage[all]{xy}
\usepackage{pstricks,graphicx}
\usepackage{stmaryrd}
\usepackage{hyperref,mdwlist,mathrsfs}
\usepackage{tikz}
\usepackage{pgf,tikz}
\usetikzlibrary{arrows}

\usepackage{color}

\usepackage{amssymb}
\usepackage{amsmath}
\usepackage{amsthm}
\usepackage{paralist}

\newcommand{\pp}{\ensuremath{\mathbb{P}}}
\newcommand{\HF}{\ensuremath{\mathrm{HF}}}

\newtheorem{theorem}{Theorem}[section]
\newtheorem{lemma}[theorem]{Lemma}

\newtheorem{proposition}[theorem]{Proposition}

\newtheorem*{proposition*}{Proposition}
\theoremstyle{definition}
\newtheorem{definition}[theorem]{Definition}
\newtheorem{example}[theorem]{Example}
\newtheorem{remark}[theorem]{Remark}

\DeclareMathOperator{\Ann}{Ann}
\DeclareMathOperator{\rank}{rank}
\DeclareMathOperator{\hess}{hess}
\DeclareMathOperator{\Hess}{Hess}

\title[Perazzo Artinian algebras]{Hilbert functions and Jordan type\\ of Perazzo Artinian algebras}

\author[N. Abdallah et al.]{Nancy Abdallah}
\address{\hspace{-15pt} Department of Engineering, Borås University, Sweden}
\email{nancy.abdallah@hb.se}

\author[]{Nasrin Altafi}
\address{\hspace{-15pt} Department of Mathematics and Statistics, Queen’s University, Kingston, Ontario, Canada and Department of Mathematics, KTH Royal Institute of Technology, Sweden}
\email{nasrinar@kth.se}

\author[]{Pietro De Poi}
\address{\hspace{-15pt} Dipartimento di Scienze Matematiche, Informatiche e Fisiche, Università degli Studi di Udine, Via delle Scienze 206, 33100 Udine, Italy}
\email{pietro.depoi@uniud.it}

 \author[]{Luca Fiorindo} 
 \address{\hspace{-15pt} Dipartimento di Matematica, Università di Genova, Via Dodecaneso 35,
16146 Genova, Italy}
  \email{luca.fiorindo@dima.unige.it}

\author[]{Anthony Iarrobino}
\address{\hspace{-15pt} Department of Mathematics, Northeastern University, Boston, MA 02115, USA}
\email{a.iarrobino@northeastern.edu}

\author[]{Pedro Macias Marques}
\address{\hspace{-15pt} Departamento de Matem\'{a}tica, Escola de Ci\^{e}ncias e Tecnologia, Centro de Investiga\c{c}\~{a}o em Matem\'{a}tica e Aplica\c{c}\~{o}es, Instituto de Investiga\c{c}\~{a}o e Forma\c{c}\~{a}o Avan\c{c}ada, Universidade de \'{E}vora, Rua Rom\~{a}o Ramalho, 59, P--7000--671 \'{E}vora, Portugal}
\email{pmm@uevora.pt}

 \author[]{Emilia Mezzetti} 
 \address{\hspace{-15pt} Dipartimento di Matematica e  Geoscienze, Universit\`a di
Trieste, Via Valerio 12/1, 34127 Trieste, Italy}
  \email{mezzette@units.it}
 
\author[]{Rosa M.\ Mir\'o-Roig} 
  \address{\hspace{-15pt} Facultat de
  Matem\`atiques i Inform\`atica, Universitat de Barcelona, Gran Via des les
  Corts Catalanes 585, 08007 Barcelona, Spain} \email{miro@ub.edu}
\author[]{Lisa Nicklasson}%
\address{\hspace{-15pt} Dipartimento di Matematica, Università di Genova, Via Dodecaneso 35,
16146 Genova, Italy}
\email{nicklasson@dima.unige.it}

\thanks{\hspace{-15pt} Altafi was supported by the grant VR2021-00472, De Poi, Fiorindo and Mezzetti are   members of INdAM - GNSAGA, Macias Marques was partially supported by FCT project UIDB/04674/2020, Mir\'o-Roig was partially supported by the grant PID2019-104844GB-I00,
Nicklasson was supported by the grant KAW-2019.0512.}

\theoremstyle{definition}



\begin{document}

\maketitle

\begin{abstract}
We study Hilbert functions, Lefschetz properties, and Jordan type of Artinian Gorenstein algebras associated to Perazzo hypersurfaces in projective space. The main focus lies on Perazzo threefolds, for which we prove that the Hilbert functions are always unimodal. Further we prove that the Hilbert function determines whether the algebra is weak Lefschetz, and we characterize those Hilbert functions for which the weak Lefschetz property holds. By example, we verify that the Hilbert functions of Perazzo fourfolds are not always unimodal.  In the particular case of Perazzo threefolds with the smallest possible Hilbert function, we give a description of the possible Jordan types for multiplication by any linear form. 
\end{abstract}

\section{Introduction}

Consider a form $F \in K[X_0, \ldots, X_N]$ of degree $d$, and its Hessian $\hess_F$, i.\,e.\ the determinant of the square matrix of second order partial derivatives of $F$. Hesse  \cite{He1}, \cite{He2} believed that the form $\hess_F$ vanishes if and only if the variety $V(F) \subset\pp^N$ is a cone, as is indeed the case when $d=2$. However, Gordan and Noether \cite{GN} proved that while Hesse's claim is true when $N \le 3$, it is false for $N>3$ and any $d \ge 3$. We refer the reader to \cite{Ru} for a review of this classical work. A class of forms $F$, discussed both in \cite{GN} and by Perazzo \cite{Pe}, with vanishing Hessian and for which $V(F)$ is not a cone, are the \emph{Perazzo forms}  $$F=X_0p_0 + X_1p_1 + \dots + X_np_n+G \in K[X_0, \ldots, X_n, U_1, \ldots, U_{m}]$$ where $n,m\geq 2$, $p_0, \ldots, p_n,G \in K[U_1, \ldots, U_{m}]$, and $p_0, \ldots, p_n$ are linearly independent but algebraically dependent.
In particular, for $N=m+n=4$, all forms with vanishing Hessian, not cones, are elements of  $K[U_1, U_2][\Delta]$ where $\Delta$ is a Perazzo polynomial of the form $X_0p_0+X_1p_1+X_2p_2$, as proved by Gordan and N\"other in \cite{GN} (see also  \cite[Theorem 7.3]{WB}).

In recent years the study of Hessians has gained new attention because of its connection to Lefschetz properties for graded Artinian Gorenstein algebras. Recall that a standard graded Artinian algebra $A$ has the \emph{weak Lefschetz property (WLP)} if multiplication by a generic linear form $\ell$ has maximal rank in each degree. Similarly $A$ has the \emph{strong Lefschetz property (SLP)} if multiplication by $\ell^s$ has maximal rank in each degree for every positive  integer $s$. Although these properties may seem easy to check, the general picture is far from understood. In the Gorenstein case, a non-trivial characterization in terms of vanishing of Hessians has been found. Indeed, the Artinian Gorenstein algebra $A_F$ with Macaulay dual generator $F$ fails the SLP if $\hess_F=0$.  More precisely, $A_F$ fails the SLP if and only if one of the  non-trivial higher Hessians of $F$ vanishes, see \cite{W} and \cite{MW}. This result has been generalised to the WLP using the so called \emph{mixed Hessians}, see \cite{GZ}.

Returning to Perazzo forms, the vanishing of $\hess_F$ implies that $A_F$ fails the SLP, and also the WLP when $d=3$. From now on, we fix ${N=4}$, which is the smallest $N$ for which a Perazzo form exists, and we tackle the problem of identifying precisely for which Perazzo polynomials $F$ the algebra $A_F$ has the WLP. If $d=4$ the WLP of $A_F$ holds for any $F$, not necessarily Perazzo, by \cite[Theorem 3.5]{Go}. The case $d \ge 5$ is studied in \cite{FMMR22} where the (componentwise) minimal and maximal possible Hilbert functions of $A_F$ are described. It is concluded that the Perazzo algebras $A_F$ with minimal Hilbert functions always satisfy the WLP, while those with maximal Hilbert function never have the WLP. For intermediate values of the Hilbert function both possibilities are realized.

In the present paper we continue the investigation of Perazzo threefolds with $d \ge 5$. We first prove in Theorem \ref{thm:unimodal} that the Hilbert function of $A_F$ is always unimodal; we then show that it always determines the presence or failure of the WLP. The unimodality result  does not generalize to higher codimension, as observed in Example \ref{ex:non-unimodal}. 
Our complete classification of the Perazzo threefolds with the WLP in terms of their Hilbert functions is given in Theorem \ref{thm:WLP}. In addition, we provide a deeper study describing the Jordan decomposition of multiplication by linear forms, both Lefschetz and non-Lefschetz, in algebras given by Perazzo threefolds with minimal Hilbert function.

The paper is organized as follows. We start by reviewing definitions and basic results concerning Hilbert functions, Lefschetz properties, Jordan types, and Perazzo algebras, in Section 2. In Section 3.1 we restrict ourselves to the case of five variables and we discover two types of Hilbert functions which force the weak Lefschetz property of Artinian Gorenstein algebras. These two results are important tools in Section 3.2, which is dedicated to the study of Hilbert functions and the WLP of Perazzo algebras. In Section 4 we discuss the possible Jordan types of multiplication by a linear form, in a Perazzo algebra with minimal Hilbert function.

\textbf{Acknowledgement.} Most of the work has been done in the \mbox{INdAM} meeting ``The Strong and Weak Lefschetz Properties'', held in\linebreak Cortona, 11--16 September 2022. The authors would like to thank the INdAM institution and the organisers for the invitation and financial support. The authors are grateful to Rodrigo Gondim and Sara Faridi for their suggestions and comments.

\section{Preliminaries}
In this section we fix notations, we recall the basic facts on Hilbert functions, Lefschetz properties, as well as on Jordan type and Perazzo hypersurfaces needed in next sections.

\subsection{Hilbert functions}
Throughout this paper  $K$ will be an algebraically closed field of characteristic zero.
Given a standard graded Artinian $K$-algebra $A=R/I$ where $R=K[x_0,x_1,\dots,x_N]$ and $I$ is a homogeneous ideal of $R$,
we denote by $\HF_A:\mathbb{Z} \longrightarrow \mathbb{Z}$ with $\HF_A(j)=\dim _KA_j=\dim _K[R/I]_j$
its Hilbert function. Since $A$ is Artinian, its Hilbert function is
captured in its \emph{$h$-vector} $h=(h_0,h_1,\dots ,h_d)$ where $h_i=\HF_A(i)>0$ and $d$ is the last index with this property. The integer $d$ is called the \emph{socle degree of} $A$.

Given  integers $n, r\ge 1$, we define the  \emph{$r$-th binomial expansion of $n$} as
\[
n=\binom{m_r}{r}+\binom{m_{r-1}}{r-1}+\cdots +\binom{m_e}{e}
\]
where $m_r>m_{r-1}>\cdots >m_e\ge e\ge 1$ 
are uniquely determined integers (see \cite[Lemma 4.2.6]{BH93}).
We write
\[ 
n^{<r>}=\binom{m_r+1}{r+1}+\binom{m_{r-1}+1}{r}+\cdots +\binom{m_e+1}{e+1}, \text{ and }
\]
\[
n_{<r>}=\binom{m_r-1}{r}+\binom{m_{r-1}-1}{r-1}+\cdots +\binom{m_e-1}{e}.
\]

The numerical functions $H:\mathbb N \longrightarrow \mathbb N$ that are Hilbert functions of graded standard $K$-algebras were characterized by Macaulay in \cite{Mac27} (see also \cite{BH93}). Indeed, given a numerical function $H:\mathbb N \longrightarrow \mathbb N$ the following conditions are equivalent:
\begin{enumerate}
\item[(i)] There exists a standard graded $K$-algebra $A$ with $H$ as Hilbert function,
\item[(ii)] $H$ satisfies the so-called \textbf{Macaulay's inequality}, i.\,e.
\end{enumerate}
\begin{equation}
\label{MacIneq}
H(0)=1\text{, and }H(t+1)\le H(t)^{<t>}\,\,\forall t\ge 1.
\end{equation}

Notice that condition (ii) imposes strong restrictions on the Hilbert function of a standard graded $K$-algebra and, in particular, it bounds its growth. Another restriction comes from the following \textbf{Green's theorem} which we recall for the sake of completeness. 

\begin{theorem} 
\label{green} 
Let $A=R/I$ be an Artinian graded algebra and let $\ell\in A_1$ be a general linear form. Let $h_t$ be the entry of degree $t$ of the h-vector of $A$. Then the degree $t$ entry $h'_t$ of the h-vector of $R/(I,\ell)$ satisfies the inequality
\[
h'_t\le (h_t)_{<t>} \text{ for all } t\ge 1.
\]
\end{theorem}
\begin{proof}
See \cite[Theorem 1]{Gr}.
\end{proof}

We recall the construction of the Artinian Gorenstein algebra $A_F$ with Macaulay dual generator a given form $F\in S=K[X_0,\dots,X_N]$; we denote by $R=K[x_0,\ldots,x_N]$ the ring of differential operators acting on the polynomial ring $S$, i.\,e.\ $x_i=\frac{\partial}{\partial X_i}$. Therefore $R$ acts on $S$ by differentiation. Given polynomials $p\in R$ and $G\in S$ we will denote by ${p\circ G}$ the differential operator $p$ applied to $G$. We define
\[
\Ann_RF:=\{p\in R \mid p\circ F=0\}\subset R,
\]
and $A_F=R/\Ann_R F$: it is a standard graded Artinian Gorenstein $K$-algebra and $F$ is called its Macaulay dual generator. We remark that every standard graded Artinian Gorenstein $K$-algebra is of the form $A_F$ for some form $F$. We may abbreviate and write ${\Ann F}$ when the ring $R$ is understood.

As an important key tool to determine Hilbert functions of Artinian Gorenstein algebras associated to Perazzo threefolds in $\pp^4$, we state the following:

\begin{proposition}
\label{exactsequence}
Let $A=A_F$ be an Artinian Gorenstein $K$-algebra and set $I=\Ann F$. Then for every linear form $\ell\in A_1$ the sequence
\begin{equation}
\label{seq}
0\longrightarrow \frac{R}{(I:\ell)}(-1)\longrightarrow \frac{R}{I}\longrightarrow \frac{R}{(I,\ell)}\longrightarrow 0
\end{equation}
is exact. Moreover $\frac{R}{(I:\ell)}$ is an Artinian Gorenstein algebra with $\ell\circ F$ as dual generator.
\end{proposition}
\begin{proof} 
We get the result cutting the exact sequence 
\[
0\longrightarrow \frac{(I:\ell)}{I}(-1)
\longrightarrow \frac{R}{I}(-1)\xrightarrow{\,\,\,\times\ell\,\,\,}\frac{R}{I}\longrightarrow \frac{R}{(I,\ell)}\longrightarrow 0
\]
into two short exact sequences. The second fact is a straightforward computation.
\end{proof}

\subsection{Lefschetz properties}

\begin{definition}
Let $A=R/I=\bigoplus_{i=0}^d\, A_i$ be a graded Artinian {$K$-algebra}. We say that $A$ has the \emph{weak Lefschetz property} (WLP, for short) if there is a linear form $\ell \in A_1$ such that, for all integers $i\ge0$, the multiplication map
\[
\times \ell: A_{i}  \longrightarrow  A_{i+1}
\]
has maximal rank, i.\,e.\ it is injective or surjective. In this case, the linear form $\ell$ is called a \emph{weak Lefschetz element} of $A$. We say that $A$ fails the WLP in degree $j$ if for a general form $\ell\in A_1$, the map $\times \ell:A_{j-1}  \longrightarrow  A_{j}$ does not have maximal rank. 

We say that $A$ has the \emph{strong Lefschetz property} (SLP, for short) if there is a linear form $\ell\in A_1$ such that, for all
integers $i\ge0$ and $k\ge 1$, the multiplication map
\[
\times \ell^k: A_{i}  \longrightarrow  A_{i+k}
\]
has maximal rank.  Such an element $\ell$ is called a \emph{strong Lefschetz element} of $A$.
\end{definition}


\begin{definition}[higher Hessians] 
Let $F\in S=K[X_0,\dots,X_N]$ be a form of degree $d$, and let $A=R/\Ann F$ be the Artinian Gorenstein $K$-algebra associated to $F$ with $h$-vector $h=(h_0,\dots,h_d)$. Given $k\le \lfloor\frac{d}{2}\rfloor$ and $\mathcal{B}=\{\alpha_i\}_i$ a $K$-basis of $A_k$, we define the \emph{$k$-th Hessian matrix} with respect to $\mathcal{B}$ to be the $h_k\times h_k$ square matrix with coefficients in $S_{d-2k}$
\[
\Hess^k_F(X_0,\dots,X_N):=\left(\alpha_i\alpha_jF(X_0,\dots,X_N)\right)_{i,j}.
\]
Moreover, we define the \emph{$k$-th Hessian} with respect to $\mathcal{B}$
\[
\hess^k_F(X_0,\dots,X_N):=\det(\Hess^k_F(X_0,\dots,X_N))\in S_{(d-2k)h_k}.
\]
\end{definition}

The definition of $k$-th Hessian and $k$-th Hessian matrix depends on the choice of a basis of $A_k$, but the vanishing of the $k$-th Hessian is independent of this choice. Thus, the condition of vanishing of the $k$-th Hessian is well-posed. The zero-th Hessian is just the polynomial $F$ and, in the case $\dim A_1=N+1$, the first Hessian, with respect to the standard basis, is the classical Hessian.
The higher Hessians are useful tools to fully characterise the strong Lefschetz elements and, in general, to understand whether an Artinian Gorenstein algebra has the strong Lefschetz property.

\begin{theorem} 
\label{thm:watanabe}
Let $F\in S$ be a homogeneous polynomial of degree $d$ and let  $ A= R/\Ann F$ be the associated Artinian Gorenstein algebra. A linear form $\ell=a_0x_0+\cdots +a_Nx_N\in A$ is a strong Lefschetz element if and only if
\[
\hess _F^k(a_0,\dots, a_N)\ne 0,\quad \text{for each} \  k=0,\dots,\left\lfloor \tfrac{d}{2}\right\rfloor.
\]        
More precisely, by fixing a $K$-basis $\mathcal{B}$ of $A_k$, since $A$ is Gorenstein, we can canonically define a basis in $A_{d-k}$. Then, up to a multiplicative constant, the $k$-th Hessian matrix $\Hess_F^k(a_0,\dots, a_N)$ is the matrix of the dual map of the multiplication map $\times \ell^{d-2k}: A_{k}  \to  A_{d-k}$.
\end{theorem}
\begin{proof}
See \cite[Theorem 4]{W} and \cite[Theorem 3.1]{MW}.
\end{proof}

\subsection{Jordan type}
\label{sec:jordantype}

\begin{definition} 
The \emph{Jordan type} $P_\ell=P_{\ell,A}$ of a graded Artinian algebra $A$  and a linear form $\ell$ of $A_1$ is the partition of $\dim_KA$ determining the Jordan block decomposition for the multiplication map $\times \ell:A\rightarrow A$. 
\end{definition}


The Jordan type of a linear form $\ell$ for  an Artinian algebra $A$ determines whether or not $\ell$ is a weak or strong Lefschetz element of $A$. In order to explain the connection we need the following definition, and for more details on Artinian algebras and Jordan type the reader is invited to look at \cite{IMM}.
\begin{definition}[Dominance order]
Given two partitions 
\[
P=(p_1,p_2,\dots, p_s), \quad p_1\ge p_2\ge \cdots \ge p_s
\]
and 
\[
Q=(q_1,q_2,\dots ,q_t), \quad q_1\ge q_2\ge \cdots \ge q_t
\]
of an integer $n$, the \emph{dominance partial order} is defined as
\[
Q\le P \hspace{2mm}\text{if} \quad \sum_{j=1}^iq_j\le \sum_{j=1}^ip_j, \hspace{2mm}\text{for all }\hspace{2mm} i\le \min\{s,t\}.
\]
\end{definition}

There is a dense open set $U\subseteq A_1$ such that if $\ell$ is a linear form in $U$ then $P_{\ell,A}\ge P_{\ell_0,A}$, for any ${\ell_0\in A_1}$. We call the Jordan type of an element in this open set the \emph{generic linear Jordan type of $A$} and denote it by $P_{A}$ (see \cite[Lemma 2.54, Definition 2.55]{IMM}; there the notation for generic linear Jordan type is $P_{1,A}$, to distinguish this from the generic Jordan type, coming from a general element of the maximal ideal).

\begin{proposition}
\label{JTandLefschetz_Prop}
Let $A$ be a graded Artinian algebra with Hilbert function $\HF_A$, then
\begin{itemize}
\item[(i)] a linear form $\ell$ is a weak Lefschetz element for $A$ if and only if the number of parts in $P_{\ell,A}$ is equal to the Sperner number of $A$; i.\,e.\ the maximum value of $\HF_A$;
\item[(ii)] a linear form $\ell$ is a strong Lefschetz element for $A$ if and only if $P_{\ell,A}$ is the conjugate partition $\HF_A^\vee$; 
\item[(iii)] for every linear form $\ell$ we have $P_{\ell,A}\le P_{A}$ w.\,r.\,t.\ the dominance order.
\end{itemize}
\end{proposition}
\begin{proof}
See \cite[Proposition 3.64]{H-W} or \cite[Proposition 2.10]{IMM}.
\end{proof}

If $P_{\ell,A}=(p_1, \ldots, p_s)$ is the Jordan type of $\ell$ then there are elements $z_1, \ldots, z_s \in A$ and \emph{strings} $(z_i, \ell z_i, \ldots, \ell^{p_i-1}z_i)$, $i=1, \ldots, s$, providing a $K$-basis for $A$. This basis is called a \emph{pre-Jordan basis} for the multiplication by $\ell$ in $A$. If in addition $\ell^{p_i}z_i=0$ for each $i$, we call it a \emph{Jordan basis}. 

\subsection{Perazzo hypersurfaces}
\label{sec:prel_perazzo}

The simplest known counterexample to Hesse's claim, i.\,e.\ a form with vanishing  Hessian which does not define a cone, is $XU^2+YUV+ZV^2$. This example was extended to a class of cubic counterexamples in all dimensions  by Perazzo in \cite{Pe}. 

\begin{definition}
\label{perazzo} 
A \emph{Perazzo hypersurface}  $X=V(F) \subset \pp^N$ is the hypersurface defined by a \emph{Perazzo form}
\[
F=X_0p_0+X_1p_1+\cdots +X_np_n+G \in K[X_0,\ldots ,X_n,U_1\ldots ,U_m]_d
\]
where $n,m\ge 2$, $N=n+m$, $p_i\in K[U_1,\ldots ,U_m]_{d-1}$ are algebraically dependent but linearly independent, and $G\in K[U_1,\ldots ,U_m]_{d}$.
\end{definition}

The fact that the $p_i$'s are algebraically dependent implies $\hess_F=0$, while the linear independence assures that $V(F)$ is not a cone. 
We are mainly interested in Perazzo hypersurfaces in $\pp^4$. In this case we use the notations 
$S=K[X,Y,Z,U,V]$ and $R=K[x,y,z,u,v]$. We have
\begin{equation}
\label{Perazzo-form}
F=Xp_0+Yp_1+Zp_2+G \quad \text{where} \ p_0, p_1, p_2, G  \in K[U,V]
\end{equation}
and any choice of $p_0, p_1, p_2$ will be algebraically dependent. The following lemma plays a key role in the induction step used in the proof of our main results of Section \ref{sec:HFforPerazzo}.

\begin{lemma}
\label{lemmapartials}
Let 
$F=Xp_0+Yp_1+Zp_2+G$ be a Perazzo form of degree $d\ge 4$  and let $A_F$ be the associated Artinian Gorenstein algebra. Then, for a general linear form $\ell\in A_F$, the polynomial $\ell\circ F$ defines a Perazzo form of degree $d-1$.
\end{lemma}
\begin{proof}
Since $V(F)$ is not a cone, we can write $\ell=a_0x+a_1y+a_2z+b_0u+b_1v$ for some coefficients $a_i,b_j\in K$ not all zero. Then we can exhibit the action of $\ell$ on $F$ as
\begin{equation*}
\begin{split}
\ell\circ F=
X\Tilde{p_0}+Y\Tilde{p_1}+Z\Tilde{p_2}+\left(a_0p_0+a_1p_1+a_2p_2+b_0\frac{\partial G}{\partial U}+b_1\frac{\partial G}{\partial V}\right)
\end{split}
\end{equation*}
with 
\begin{gather*}
\Tilde{p_0}=b_0\frac{\partial p_0}{\partial U}+b_1\frac{\partial p_0}{\partial V},\quad \Tilde{p_1}=b_0\frac{\partial p_1}{\partial U}+b_1\frac{\partial p_1}{\partial V},\\ 
\text{ and } \quad
\Tilde{p_2}=b_0\frac{\partial p_2}{\partial U}+b_1\frac{\partial p_2}{\partial V}.
\end{gather*}

The form $\ell\circ F$ has degree $d-1\ge 3$. It remains to prove that the polynomials $\Tilde{p_0},\Tilde{p_1},\Tilde{p_2}$ are linearly independent. Since $\ell\circ F$ corresponds to an adjoint surface of $V(F)$ which is not a cone, then, for the general linear form, $\ell\circ F$ does not define a cone \footnote{This may not happen for some specific linear form (e.g. $\ell=x+y+z$).}.
\end{proof}

\section{Hilbert functions occurring for Perazzo algebras, and the weak Lefschetz property}\label{sec:HFforPerazzo}
This section is devoted to the proof of our first main result Theorem \ref{thm:WLP}. We first prove some general results about the WLP of Artinian Gorenstein algebras. 

\subsection{The weak Lefschetz properties for some Artinian Gorenstein algebras}

Let $A=R/I$ be a graded Artinian Gorenstein algebra. In the following two propositions we prove that specific conditions on the Hilbert function of $A$ ensure the WLP.

\begin{proposition}
\label{prop:gor_linearHF_WLP}
Let $A=R/I$ be an Artinian Gorenstein algebra of even socle degree $d=2s$ with Hilbert function
 $$H=(1,h_1,h_2,\dots,h_{s-1},h_s,h_{s-1},\dots, h_2,h_1,1)$$ with $h_1,\dots,h_s$ consecutive increasing integers and $h_k\leq k$ for some $\frac{d}{2}+1\leq k\leq d-1$. Then $A$ has the WLP.
\end{proposition}

\begin{proof} 
Consider the exact sequence (\ref{seq}) in Proposition~\ref{exactsequence} where $\ell$ is a general linear form. Denote by $h'_k=\dim [R/(I,\ell)]_k$ and $\tilde{h}_k=\dim [R/(I:\ell)]_k$. Then $h_k-\Tilde{h}_{k-1}=h'_k$ for all $k>0$ and $h_0=h'_0=1$. Let $k_0$ be the smallest $k$ such that $h_k\leq k$. By Theorem~\ref{green} $h'_k=0$ for all $k\ge k_0$, and therefore, $\Tilde{h}_{k}=h_k$ for all $k\ge k_0$. 
Since $R/I$ and $R/(I:\ell)$ are Gorenstein, by symmetry we have $h'_k=h_k-\Tilde{h}_{k-1}=h_k-(h_{k}-1)=1$ for all $k\leq d-k_0+1$. Hence, Macaulay's inequality (\ref{MacIneq}) forces  $h'_{s+1}\leq 1$. If $(h'_{s-1},h'_s,h'_{s+1})=(1,1,1)$ then $\bigl(\tilde{h}_{s-2},\tilde{h}_{s-1},\tilde{h}_{s}\bigr)=(1,2,1)$ which again contradicts both Macaulay's theorem and the Gorenstein property of $R/(I:\ell)$. Hence, $h'_{s+1}=0$ and $R/I$ has the WLP. 
\end{proof}

\begin{example}\label{ex:WLP_h-vec}
Every Artinian Gorenstein algebra with $h$-vector
$$(1,5,6,7,\dots, s-1,s,s-1, \dots, 7,6,5,1)$$ 
and socle degree at least 5 has the WLP.
\end{example}
As another example of Hilbert function forcing WLP we have:

\begin{proposition}
\label{prop:flat_top_WLP}
Every Artinian Gorenstein algebra with $h$-vector 
$$ h_0<h_1<\cdots<h_{t-1}<h_t=\cdots=h_s>h_{s+1}>\cdots>h_{d-1}>h_d$$
where $s\ge t+2$ and $h_s\le s$ has the WLP.
 \end{proposition}
 
\begin{proof}
It is enough to check that for a general $ \ell\in R_1$  the multiplication map $\times \ell: A_{d/2} \rightarrow  A_{d/2+1}$ is surjective if $d$ is even, resp.\ 
$\times \ell:A_{(d-1)/2} \rightarrow A_{(d+1)/2}$ if $d$ is odd. So, it suffices to prove that $[R/(I,\ell)]_t =0$ for $t\ge (d+1)/2$.

Since $h_s\le s$, applying Green's theorem we get
$[R/(I,\ell)]_s=0$ i.\,e.\ the multiplication map $\times \ell:A_{s-1} \rightarrow  A_{s}$ is an isomorphism. By duality
the multiplication map $\times \ell:A_{t} \rightarrow A_{t+1}$ is also an isomorphism. This implies that $[R/(I,\ell)]_{t+1}=0$. Therefore,
$[R/(I,\ell)]_i=0$ for $i\ge t+1$ and we are done.
\end{proof}

\begin{remark}
For a graded Artinian Gorenstein algebra $A$ of codimension $n$ and $h$-vector satisfying $h_t=h_{t+1}\le t$, for some $t$, a result of Gotzmann shows that $A$ is the quotient of the coordinate ring of a zero-dimensional scheme in $\mathbb{P}^{n-1}$ of degree $h_t$, see \cite{G} or \cite[Proposition C.32]{IK} for more details. Therefore, such an algebra $A$ has the WLP. 

The special case of Proposition \ref{prop:flat_top_WLP} when $h_{s-1}=h_s\le s-1$ follows from the result of Gotzmann. 
\end{remark}
\begin{example} Every Artinian Gorenstein algebra with $h$-vector $$(1, 4, 6, 7^k, 6, 4, 1), \quad \text{for} \  k\ge5$$ has the WLP. Gotzmann's result proves it for every $k\geq 6$ and the case $k=5$ is implied by Proposition \ref{prop:flat_top_WLP}.
\end{example}

\subsection{Hilbert functions and WLP for algebras associated to a Perazzo hypersurface }

We now continue with the case $A_F$ being an Artinian Gorenstein algebra associated to a Perazzo hypersurface of degree $d \ge 5$ in $\pp^4$. Recall that by \cite[Theorem 4.3]{FMMR22} the algebra $A_F$ has the weak Lefschetz property if the Hilbert function of $A_F$ is the termwise minimal one, namely $(1,5,6, \ldots, 6,5,1)$. We also know that $A_F$ fails the WLP if the Hilbert function is the termwise maximal for the given socle degree, by \cite[Theorem 4.1]{FMMR22}. We will now extend the two results and give a complete classification of such algebras $A_F$ with the WLP. The case $d=5$ is completely covered by the results of \cite{FMMR22}, as the only possible Hilbert functions are the minimal, $(1,5,6,6,5,1)$, and the maximal, $(1,5,7,7,5,1)$. Moreover, we have the following general result. 

\begin{theorem}
\label{thm:unimodal} 
The Hilbert function of an Artinian Gorenstein algebra associated to a Perazzo hypersurface of degree $d\ge 5$ in $\pp ^4$ is unimodal.
\end{theorem}
\begin{proof} 
By \cite[Proposition 3.5]{FMMR22} we have $h_i<(i+3)(2d-3i)/2$ for $i\le (d/2) - 1$. Indeed, for $1\le i\le (d/2) - 1$ we have
\[
h_i\le 4i+1<\tfrac{1}{2}(i+3)(8i-2-3i)\le\tfrac{1}{2}(i+3)(2d-3i),\quad\text{ for }i\le \tfrac{d+1}{4}.
\]
\[
h_i\le d+2<\tfrac{1}{2}(i+3)(4i+4-3i)\le\tfrac{1}{2}(i+3)(2d-3i),\quad\text{ for }i > \tfrac{d+1}{4}.
\]
Applying \cite[Proposition 2.6]{MNZ} we conclude that $h_{i+1}\ge h_i$ and hence we have the unimodality of the Hilbert function of $A_F$. 
\end{proof}

Theorem \ref{thm:unimodal} does not generalize to $\pp^N$ with $N >4$.
\begin{example}\label{ex:non-unimodal}
Let $F=Xp_0+Yp_1+Zp_2$ where $p_0,p_1,p_2$ are general forms of degree $9$ in the variables $U,V,W$. Then the algebra $A_F$ has Hilbert function 
\[
(1, \ 6, \ 15, \ 28, \ 43, \ 42, \ 43, \ 28, \ 15, \ 6, \ 1) 
\]
which is not unimodal. 
\end{example}

In the proof of Proposition \ref{prop:WLP_d6} we need to exclude certain vectors as possible $h$-vectors of Perazzo algebras. It was pointed out by Mats Boij that the Hilbert function $(1,5,6,8,6,5,1)$ cannot occur for any Artinian Gorenstein algebra, and this can be proven by using Betti tables of lexicographic ideals and a cancellation argument. We state the result below. 

\begin{lemma}\label{lemma:excluded_hvector}
There is no Artinian Gorenstein algebra with Hilbert function $(1,5,6,8,6,5,1)$.
\end{lemma}

\begin{proof} 
Let $A=R/I$ be an Artinian Gorenstein algebra with Hilbert function $(1,5,6,8,6,5,1)$. In \cite{Mac27}, Macaulay proved that there exists a lexicographic ideal $I_{\operatorname{lex}}$ with the same Hilbert function.  By \cite[Theorem 1.1]{P}, the Betti numbers $\beta_{ij}(S/I)$ must arise from $\beta_{ij}(S/I_{\operatorname{lex}})$ by consecutive cancellations. The minimal Betti table of $S/I_{\operatorname{lex}}$ is 
{\scriptsize
\begin{verbatim}
				+-----------------------+
				|       0  1  2  3  4  5|
				|total: 1 22 67 84 49 11|
				|    0: 1  .  .  .  .  .|
				|    1: .  9 20 20 10  2|
				|    2: .  2  5  4  1  .|
				|    3: .  4 15 21 13  3|
				|    4: .  2  7  9  5  1|
				|    5: .  4 16 24 16  4|
				|    6: .  1  4  6  4  1|
				+-----------------------+
\end{verbatim}%
}%
\noindent
Here we use the convention that the entry at row $i$ and column $j$ is the Betti number $\beta_{j,i+j}$.
   To get an Artinian Gorenstein algebra $A$, the Betti number ${\beta_{56}(S/I_{\operatorname{lex}})=2}$ should be cancelled to zero. This can only happen if $\beta_{46}(S/I_{\operatorname{lex}})$ or $\beta_{66}(S/I_{\operatorname{lex}})$ is at least $2$, a contradiction. 
\end{proof}

The following proposition serves as the basis for an induction in the proof of our general Theorem \ref{thm:WLP}.

\begin{proposition}\label{prop:WLP_d6}
Let $A_F$ be an Artinian Gorenstein algebra of socle degree 6 associated to a Perazzo hypersurface $V(F)\subset \pp^4$. The algebra $A_F$ fails WLP if and only if the Hilbert function is maximal, i.\,e. $(1, 5, 8, 8, 8, 5, 1)$.
\end{proposition}

\begin{proof} 
By Theorem \ref{thm:unimodal} and Lemma \ref{lemma:excluded_hvector}
\begin{gather*}
    (1,5,6,6,6,5,1),\quad (1,5,6,7,6,5,1),\quad (1,5,7,7,7,5,1),\\ 
    \quad
    (1,5,7,8,7,5,1), \quad\text{and}\quad (1,5,8,8,8,5,1).
\end{gather*}
are the possible $h$-vectors for  $A_F$.
We know that $A_F$ has the WLP if its Hilbert function is the minimal $(1,5, 6, 6, 6,5,1)$.
The case of $(1,5,6,7,6,5,1)$ is covered by Example \ref{ex:WLP_h-vec}. We also know that $A_F$ fails the WLP if it has the maximal Hilbert function $(1,5,8,8,8,5,1)$. 
Let us now see that the $h$-vector $(1,5,7,7,7,5,1)$ can be excluded.


Assume that $A_F=R/I$ has the Hilbert function $(1,5,7,7,7,5,1)$.  Take $\ell $ a general linear form and denote by $h'$ the Hilbert function  of $R/(I,\ell)$. As in the previous case, the Hilbert function of $R/(I:\ell)$ is $(1,5,a,a,5,1)$ with $a=6$ or $7$. We consider another general linear form $\ell'\in A_1$ and the commutative diagram
\[
\begin{array}{ccccccccc}
0 & \rightarrow & [R/(I:\ell)]_2 & \rightarrow & A_3=[R/I]_3 & \rightarrow & [R/(I,\ell)]_3 & \rightarrow & 0\\ 
& & \downarrow & & \downarrow & & \downarrow \\
0 & \rightarrow & [R/(I:\ell)]_3 & \rightarrow & A_4=[R/I]_4 & \rightarrow & [R/(I,\ell)]_4 & \rightarrow & 0
\end{array}
\]
where the vertical morphisms are given by the multiplication by $\ell'$. If $a=6$, $R/(I:\ell)$ has minimal Hilbert function $(1, 5, 6, 6, 5, 1)$. Therefore, $R/(I:\ell)$ has the WLP and the first vertical map is an isomorphism. $R/(I,\ell)$ has Hilbert function $(1, 4, 2, 1, 1)$ and so the last vertical map is surjective while the middle map is not surjective because $A_F$ fails the WLP in degree $4$.
If $a=7$, $R/(I:\ell)$ has maximal Hilbert function $(1, 5, 7, 7, 5, 1)$. Therefore, $R/(I:\ell)$ fails the WLP while $A_F$ has the WLP.  Hence the middle map is an isomorphism which contradicts the fact that the first vertical map is not injective.

It only remains to prove WLP in the case that $A_F$ has the Hilbert function $(1,5,7,8,7,5,1)$.
In that case, let $\ell$ be a general linear form, and note that $A_F$ has the WLP if and only if $\times \ell: [A_F]_3 \to [A_F]_4$ is surjective. 
 
Write  $A_F=R/I$  and consider the algebra $R/(I:\ell)$. By Lemma \ref{lemmapartials} this is a Gorenstein algebra associated to a Perazzo form of degree five. Recall that the Hilbert function of $R/(I,\ell)$ is determined by the Hilbert functions of $R/I$ and $R/(I:\ell)$ as explained in Proposition \ref{exactsequence}. If $R/(I:\ell)$ has the Hilbert function $(1,5,7,7,5,1)$ then $R/(I,\ell)$ has the Hilbert function $(1,4,2,1)$. In particular $[R/(I,\ell)]_4=0$, and we are done. The only other possibility for the Hilbert function of $R/(I:\ell)$ is the minimal $(1,5,6,6,5,1)$, in which case $R/(I,\ell)$ would have the Hilbert function $(1,4,2,2,1)$. 
We consider another general linear form $\ell'\in A_1$ and the commutative diagram
\[
\begin{array}{ccccccccc}
0 & \rightarrow & [R/(I:\ell)]_2 & \rightarrow & A_3=[R/I]_3 & \rightarrow & [R/(I,\ell)]_3 & \rightarrow & 0\\ 
& & \downarrow & & \downarrow & & \downarrow \\
0 & \rightarrow & [R/(I:\ell)]_3 & \rightarrow & A_4=[R/I]_4 & \rightarrow & [R/(I,\ell)]_4 & \rightarrow & 0
\end{array}
\]
where the vertical morphisms are given by the multiplication by $\ell'$. The algebra $R/(I:\ell)$ has the WLP and hence
\[
\times \ell':[R/(I:\ell)]_2 \longrightarrow [R/(I:\ell)]_3
\]
is an isomorphism. Moreover, it is easy to see that 
\[
\times \ell':[R/(I,\ell)]_3\longrightarrow [R/(I,\ell)]_4
\]
is surjective as $[R/(I,\ell)]_4$ is a one dimensional space. By the snake lemma
\[
\times \ell':[R/I]_3\longrightarrow [R/I]_4
\]
is surjective, which implies the WLP. 
\end{proof}

\begin{remark}
    Note that the Hilbert function $(1,5,7,7,7,5,1)$ cannot occur for a Perazzo hypersurface, but can still occur for some Artinian Gorenstein algebra. Indeed, it is the Hilbert function of the algebra defined by the Macaulay dual generator $F = X_0^6+X_0X_1^5+X_0X_2^5+X_1^6+X_2^6+U^6+V^6.$
\end{remark}

We can now give the classification of Artinian Gorenstein algebras associated to Perazzo hypersurfaces of degree $d \ge 5$ in $\pp^4$ with the weak Lefschetz property. Recall that such algebras have Sperner number at most $d+2$.
\begin{theorem}\label{thm:WLP} 
Let $A_F$ be an Artinian Gorenstein algebra associated to a Perazzo hypersurface $V(F)\subset \pp^4$ of degree $d\ge 5$. Let $(h_0,h_1,\ldots ,h_d)$ be its Hilbert vector. The algebra $A_F$ has the WLP if and only if $\# \{i \mid h_i=d+2 \} \le 1$.
\end{theorem}
\begin{proof}
Assume first that $\# \{i \mid h_i=d+2 \} \ge 2$ (resp. $\ge 3$) if $d$ is odd (resp. even). 
Consider the Hessian of $F$ of order $(d-1)/2$ (resp. $(d-2)/2)$. To compute it we fix a monomial basis $\mathcal B$ of the homogeneous component of $A_F$ of degree $(d-1)/2$ (resp. $(d-2)/2$), formed by $d+2$ monomials. Since there are at most $(d+1)/2$ monomials involving only $u,v$ (resp. $d/2$), the number of monomials containing some of the variables $x, y, z$ is at least $(d+3)/2$ (resp. $(d+4)/2$). This implies that in the Hessian matrix there is a block of zeros of order at least $(d+3)/2$ (resp. $(d+4)/2$), which forces the determinant to vanish. Therefore the multiplication map by a general linear form does not have maximal rank in degree $(d+1)/2$ (resp. in degree $d/2$), and $A_F$ fails the WLP.
 
Let us now assume $\# \{i \mid h_i=d+2 \} \le 1$, and prove that $A_F$ has the WLP. We have seen that the result holds for $d=5,6$. We proceed by induction over $d$, treating the cases of $d$ even and $d$ odd separately.

 \vskip 2mm

\noindent 
\textbf{CASE $d$ even}. Write $d=2s$. The hypothesis $\# \{i \mid h_i=d+2 \} \le 1$ together with Theorem \ref{thm:unimodal} implies that $h_{s+1}\le d+1$. We let $I=\Ann F$ so that $A_F=R/I$. 
Take a general linear form $\ell\in A_F$ and consider the short exact sequence (\ref{seq}).
By Lemma \ref{lemmapartials} we know that $R/(I:\ell)$
is an Artinian Gorenstein algebra associated to the Perazzo form $\ell\circ F$ of degree $d-1=2s-1$. Denote by $h_i'$ and $\widetilde{h_i}$ the Hilbert functions of $R/(I,\ell)$ and $R/(I:\ell)$, and recall that $\widetilde{h_i}=h_{i+1}-h'_{i+1}$. Using Green's theorem (see Theorem \ref{green}) and the inequality $h_{s+1}<d+2$, we get that $h'_{s+1} \le 1$. If $h'_{s+1}=0$ the map 
$$\times \ell: [R/I]_s\longrightarrow [R/I]_{s+1}$$
is surjective, and $A_F$ has the WLP. Suppose instead $h'_{s+1}=1$. Then, as $h_{s+1}\le d+1$, we get $\widetilde{h_{s}}<d+1$. Therefore, we can apply our induction hypothesis to 
$R/(I:\ell)$ and conclude that $R/(I:\ell)$ satisfies the WLP. Consider now another general linear form $\ell'\in R/I$ and the commutative diagram
\begin{equation}
\label{commutativediagram}
\begin{array}{ccccccccc}
0 & \rightarrow & [R/(I:\ell)]_{i-1} & \rightarrow & [R/I]_i & \rightarrow & [R/(I,\ell)]_i & \rightarrow & 0\\
& & \downarrow & & \downarrow & & \downarrow\\
0 & \rightarrow & [R/(I:\ell)]_i & \rightarrow & [R/I]_{i+1} & \rightarrow & [R/(I,\ell)]_{i+1} & \rightarrow & 0
\end{array}
\end{equation}
for $i=s$, and where the vertical morphisms are given by the multiplication by $\ell'$. Since $R/(I:\ell)$ satisfies the WLP, the leftmost vertical map is an isomorphism.
On the other hand it is easy to check, as $h_{s+1}'=1$, that the rightmost vertical map is an epimorphism. By the snake lemma, the middle vertical map is an epimorphism, which then gives $A_F$ the WLP.

\vskip 2mm
\noindent 
\textbf{CASE $d$ odd}. Write $d=2s+1$. By Theorem \ref{thm:unimodal} we know that $h$ is symmetric and unimodal. Since we assume $\# \{i \mid h_i=d+2 \} \le 1$ we have $h_s=h_{s+1}\le d+1$ in this case.  
Applying Theorem \ref{green} we get $h'_{s+1} \le 1$. If $h'_{s+1}=0$ we are already done, so suppose $h'_{s+1} = 1$. By Macaulay's inequality then $h'_i \le 1$ for all $i>s$. Considering the commutative diagram \eqref{commutativediagram} we see that the rightmost vertical map is surjective for $i>s$. Since $R/(I:\ell)$ is a Perazzo of socle degree $2s$ it has the WLP, and therefore the leftmost vertical map in \eqref{commutativediagram} is surjective for $i>s$. We can conclude that the middle map $\times \ell':[R/I]_i \to [R/I]_{i+1}$ is surjective when $i>s$, and by duality also injective for $i<s$. This gives the equality $h'_{i}=h_{i}-h_{i-1}$ for $i \le s$. 
 
Next, we claim that $h_i'=1$ for some $i \le s$, or equivalently $h_i=h_{i-1}+1$. To prove the claim, consider the possibility that $h_i\ge h_{i-1}+2$ for each $i=2, \ldots, s$. Then we would have $h_s \ge h_1+2(s-1)=5+2(s-1)=d+2$ which would contradict the assumption $h_s \le d+1$. 
 
Now Macaulay's inequality together with the fact that $h'_i=1$ for some $i \le s$ implies $h'_s=1$. Finally, we consider the diagram \eqref{commutativediagram} with $i=s$. In this case the leftmost vertical map is injective, and the rightmost is bijective as $h_s'=h'_{s+1}=1$. We conclude that the middle map is injective, and therefore $R/I$ has the WLP. 
\end{proof}
 \begin{theorem}
\label{unimodalhyperplanesection} 
Let $A_F$ be an Artinian Gorenstein algebra associated to a Perazzo hypersurface $V(F)\subset \pp ^4$ of degree $d\ge 5$. Assume that $A_F$ has the WLP and that $\ell$ is a weak Lefschetz element. Then the  Hilbert function of $A_F/(\ell )$ is unimodal.
\end{theorem}
\begin{proof}
Denote the $h$-vectors of $A_F$ and $A_F/(\ell )$ by $h$ and $h'$, respectively. As $\ell$ is a weak Lefschetz element we have $h' = \Delta h$, meaning $h'_i=\max(h_i-h_{i-1},0)$. Using \cite[Proposition 3.4]{FMMR22} we get that $h'\le 4$ for every $i$. Macaulay's inequality (\ref{MacIneq}) on $h'$, for every $i\ge 2$, implies that $(h')_{i+1}\leq h'_i$ unless $i=2$, $h'_2=3$ and in that case $h'_3\leq 4$. This means that $h'$ is unimodal except when $h' = (1,4,3,4,\ldots)$. We will show that there is no Artinian Gorenstein algebra $A_F$ associated to a Perazzo hypersurface $V(f)\subset \pp ^4$ such that $\Delta h = (1,4,3,4,\ldots)$. If such an algebra would exist, it would have Sperner number at least 12, which implies $d \ge 10$. But for $d=10$ the vector $h$ does not satisfy the necessary condition for WLP given in Theorem \ref{thm:WLP}. Suppose $d \ge 11$. Then $h = (1,5,8,12,h_4,h_5,\dots , h_5,h_4,12,8,5,1)$. We will show that $h$ cannot occur as the Hilbert function of $A_F$ by showing that $h_3=12$ forces $h_2=9$. 

We  denote the  polynomial $F$ as
$$F = Xp_0(U,V)+Yp_1(U,V)+Zp_2(U,V)+G(U,V),$$ where 
\begin{align*}
p_0(U,V) = \sum_{i=0}^{d-1}\tbinom{d-1}{i}a_iU^{d-1-i}V^i, &\quad  p_1(U,V) = \sum_{i=0}^{d-1}\tbinom{d-1}{i}b_iU^{d-1-i}V^i,\\
p_2(U,V) = \sum_{i=0}^{d-1}\tbinom{d-1}{i}c_iU^{d-1-i}V^i,&\quad \text{and}\quad 
G(U,V) = \sum_{i=0}^{d}\tbinom{d}{i}g_iU^{d-i}V^i.
\end{align*}

From \cite[Proposition 3.4]{FMMR22} we obtain that  
\begin{align}
        &\rank  M_2 + \rank  N_2\le h_2 \le \rank  M_2 + \rank  N'_2,\hspace{2mm}\text{and} \label{h2Eq} \\ 
        &\rank  M_3 + \rank  N_3\le h_3 \le \rank  M_3 + \rank  N'_3, \label{h3Eq} 
\end{align}

\noindent where $M_2, M_3,N_2, N_3, N^\prime_2,$ and $N^\prime_3$ are the matrices
\begin{align*}
&M_2 = \left( A_2 \vert B_2\vert C_2\right), \quad M_3 = \left( A_3 \vert B_3\vert C_3\right),\\
&{N_2=\begin{pmatrix} A_{3} \\ \hline B_{3} \\ \hline C_{3} \end{pmatrix}, \quad N_3 = \begin{pmatrix} A_{4} \\ \hline B_{4} \\ \hline C_{4} \end{pmatrix}}, \quad N^\prime_2=\begin{pmatrix} A_{3} \\ \hline B_{3} \\ \hline C_{3} \\ \hline G_2 \end{pmatrix}, \quad \text{and} \quad N^\prime_3 = \begin{pmatrix} A_{4} \\ \hline B_{4} \\ \hline C_{4} \\ \hline G_3 \end{pmatrix}
\end{align*}
with building blocks
\begin{Small}
\begin{align*}
A_k = \begin{pmatrix} a_0&a_1&\dots &a_{k-1}\\
a_1&a_2&\dots &a_{k-2}\\
\vdots&\vdots &\cdots &\vdots \\
a_{d-k}&a_{d-k+1}&\dots &a_{d-1}\\
\end{pmatrix},\quad & B_k = \begin{pmatrix} b_0&b_1&\dots &b_{k-1}\\
b_1&b_2&\dots &b_{k-2}\\
\vdots&\vdots &\cdots &\vdots \\
b_{d-k}&b_{d-k+1}&\dots &b_{d-1}\\
\end{pmatrix},\\
C_k = \begin{pmatrix} c_0&c_1&\dots &c_{k-1}\\
c_1&c_2&\dots &c_{k-2}\\
\vdots&\vdots &\cdots &\vdots \\
c_{d-k}&c_{d-k+1}&\dots &c_{d-1}\\
\end{pmatrix},\quad 
 & G_k = \begin{pmatrix} g_0&g_1&\dots &g_{k}\\
g_1&g_2&\dots &g_{k-1}\\
\vdots&\vdots &\cdots &\vdots \\
g_{d-k}&g_{d-k+1}&\dots &g_{d}\\
\end{pmatrix}.
\end{align*}
\end{Small}
We claim that $\rank M_3\ge 8$ implies that $\rank M_2 = 6$. To show the claim notice that the matrix obtained by removing the first row (or the last row) from $M_2$ is a submatrix of $M_3$. Suppose $\rank M_2< 6$ and assume that the last column of $M_2$ is 
a linear combination of the other columns. This implies that the last two columns of $M_3$ are in the span of the first $7$ columns. So $\rank M_3\leq 7$ and this proves the claim. \par
Observe that $\rank N^\prime_3\le 4$, so $h_3=12$ implies that $\rank M_3 \ge 8$ by inequality (\ref{h3Eq}). In this case the claim shows that $M_2$ has to have maximal rank, $\rank M_2=6$. On the other hand, $\rank M_3\ge 8$ forces at least one of the three blocks $A_3, B_3,$ or $C_3$ to have maximal rank, and therefore $N^\prime_2$ and $N_2$ both have maximal rank, that is equal to three. So $\rank M_2+\rank N_2 = \rank M_2+\rank N^\prime_2=6+3=9$ and by inequality (\ref{h2Eq}) we conclude that $h_2=9$ which completes the proof.
\end{proof} 

Theorem \ref{thm:WLP} and Theorem \ref{unimodalhyperplanesection} rule out some intermediate possibilities for $h$-vectors of Perazzo algebras in between the maximum and minimum $h$-vectors given by \cite[Propositions 3.5 and 3.6]{FMMR22}. 
In general, it is a difficult problem to determine whether a $h$-vector is the Hilbert function of an Artinian Gorenstein algebra, let alone one coming from a Perazzo hypersurface.

 \begin{example}\label{ex:imposHF_d6}
There is no Artinian Gorenstein algebra associated to a Perazzo hypersurface in $\pp^4$ with Hilbert function $$H=(1,5,6,7,9,9,9,7,6,5,1).$$ 
If there were an algebra with this Hilbert function, then Theorem \ref{thm:WLP} would imply that it has WLP, and by Theorem \ref{unimodalhyperplanesection} $\Delta H$ should  be unimodal, which is a contradiction.
\end{example}

\section{Jordan type of a Perazzo algebra}
\label{JTPerazzo}

In this section, we study the Jordan types of an Artinian Gorenstein algebra corresponding to a Perazzo threefold of minimal Hilbert vector, i.\,e.\ of the type $(1,5,6,\dots,6,5,1)$. 

\subsection{Perazzo algebras with minimal Hilbert function}
An explicit classification of the possible dual generators $F$ of degree $d \ge 5$ defining a Perazzo threefold with minimal Hilbert function is given in \cite[Theorem 5.4]{FMMR22}. We state a slightly rephrased version of this result. 

\begin{lemma}\label{lemma:min_perazzo}
    Let $F \in K[X,Y,Z,U,V]$ be a Perazzo form such that the algebra $A_F$ has minimal Hilbert function $(1,5,6, 6, \ldots, 6,5,1)$. Then the dual generator $F$ can be expressed as
    \begin{enumerate}[(i)]
    \item \label{minimalHFPerazzo1}$XU^{d-1}+YU^{d-2}V+ZU^{d-3}V^2$, 
    \item \label{minimalHFPerazzo2}$XU^{d-1}+YU^{d-2}V+ZV^{d-1}$, or
    \item \label{minimalHFPerazzo3} $XU^{d-1}+Y(U+\lambda V)^{d-1}+ZV^{d-1}$ with $\lambda\in K^*$
\end{enumerate}
after a linear change of variables. 
\end{lemma}
\begin{proof}
    By \cite[Theorem 5.4]{FMMR22}, there are three classes of forms up to a linear change of variables. In the first case, $F$ can be written as
    \begin{align*}
        F&=XU^{d-1}+YU^{d-2}V+ZU^{d-3}V^2+aU^d+bU^{d-1}V+cU^{d-2}V\\
        &=(X+aU)U^{d-1}+(Y+bU)U^{d-2}V+(Z+cU)U^{d-3}V^2
    \end{align*}
    for $a,b,c\in K$. A further linear change of variables gives the case (i). The other cases are obtained in the same manner.
\end{proof}
 
We know by \cite[Theorem 4.3]{FMMR22} that, in all the three cases of Lemma \ref{lemma:min_perazzo}, $A_F$ has the WLP. Following the proof in \cite[Theorem 3.25]{Luca} we get the precise description of the weak Lefschetz elements. 

\begin{proposition}
\label{prop:Lefschetz-elements} 
For Perazzo threefolds with minimal Hilbert function, the Lefschetz elements are the linear forms $a_0x+a_1y+a_2z+b_0u+b_1v$ satisfying the following conditions, for the three types of $F$ listed in Lemma \ref{lemma:min_perazzo}:
\begin{enumerate}[(i)]
    \item  $b_0\ne 0$,
    \item $ b_0b_1 \ne 0$,
    \item $ b_0b_1(b_0+\lambda b_1) \ne 0 .$
\end{enumerate}
\end{proposition}

In our situation, the Jordan type of the algebra $A_F$ is a partition of $6d-6$.  Let us start with the generic linear Jordan type $P_{A_F}=(p_1, \ldots, p_t)$. We know that a general linear form $\ell\in A_F$ satisfies $\ell^d\ne0$, so ${p_1=d+1}$. 
By Proposition \ref{prop:Lefschetz-elements} above, a general linear form $\ell \in A_F$ is a weak Lefschetz element.
So, since the Sperner number is $6$, this must be the number of parts in $P_{A_F}$ (see Proposition \ref{JTandLefschetz_Prop}). Since $A_F$ does not have the SLP, $P_{A_F}$ is strictly dominated by the conjugate of the Hilbert function:  $H(A_F)^\vee=\bigl(d+1,(d-1)^4,d-3\bigr)$. The only possibility is  
\begin{equation}
\label{genericlinearJTeq}
P_{A_F}=\bigl(d+1,(d-1)^3,(d-2)^2\bigr),
\end{equation}
as any other partition of $6d-6$ with $p_1=d+1$ and strictly dominated by $H(A_F)^\vee$ has more than $6$ parts. 

\begin{remark}
If we consider a general element $\ell'$ in the maximal ideal of $A_F$ (not necessarily homogeneous), we could ask what its Jordan type is (this is the generic Jordan type of $A_F$, see \cite[Definitions 2.1 and 2.55]{IMM}). We know that this Jordan type is dominated by $H(A_F)^\vee$ and since a general element in the maximal ideal of $A_F$ specialises to a general linear form in $(A_F)_1$, its Jordan type dominates the generic linear Jordan type (see the discussion before Lemma 2.54 in \cite{IMM}).
But if we compare the partition $H(A_F)^\vee$ with the generic linear Jordan type we have just obtained, we see that there is no partition between the two in the dominance order, so the generic Jordan type must be equal to one of these partitions. Since the Jordan type $H(A_F)^\vee$ implies that the vector space $(\ell')^{d-1}(A_F)_1$ has dimension $5$ and this can only be attained if the same happens for $(\ell'_1)^{d-1}(A_F)_1$, where $\ell'_1$ is the linear summand of $\ell'$, we see that the generic Jordan type of $A_F$ is $P_{A_F}$ as in \eqref{genericlinearJTeq}. 
\end{remark}

Let us now consider the possible Jordan types for multiplication by any linear form in $A_F$.

\begin{theorem}
\label{thm:JT}
Let $A_F$ be an Artinian Gorenstein algebra of a Perazzo threefold with minimal Hilbert function. 
The two Jordan types 
\[
\bigl(d+1,(d-1)^3,(d-2)^2\bigr)\quad \text{and} \quad 
\bigl(d^2,(d-1)^2,(d-2)^2\bigr) 
\]
occur for multiplication by a weak Lefschetz element. For elements in the non-Lefschetz locus, the following holds regarding the Jordan basis, for the three types of $F$ given in Lemma \ref{lemma:min_perazzo}.
\begin{enumerate}[(i)]
    \item All strings are of lengths at most $4$,
    \item The strings have lengths $\ge d-1$ or $\le 3$,
    \item The strings have lengths $\ge d-1$ or $\le 2$.
\end{enumerate}
\end{theorem}
\begin{proof}
We have already seen that the Jordan type of a general linear form is $\bigl(d+1,(d-1)^3,(d-2)^2\bigr)$, and that this is the only possible Jordan type for a weak Lefschetz element $\ell$ with $\ell^d \ne 0$. 
Let's consider a weak Lefschetz element $\ell$ such that $\ell^d=0$. For instance, take $\ell=b_0u+b_1v$ with $b_0$, $b_1$ satisfying the conditions given in Proposition \ref{prop:Lefschetz-elements}. In this case the Jordan strings have length at most $d$. We claim that the largest part of such a Jordan type is equal to $d$. By Proposition \ref{JTandLefschetz_Prop} every Jordan type partition has to be smaller than the generic Jordan type $\bigl(d+1,(d-1)^3,(d-2)^2\bigr)$ w.\,r.\,t.\ the dominance order, so there is no partition of $6d-6$ with exactly $6$ parts that is dominated by the generic Jordan type and has parts of length at most $d-1$. Thus the largest part has length equal to $d$ and by symmetry \cite[Lemma 4.6]{CG}
we must have at least $2$ parts of length $d$. So, the Jordan type in this case is of the form $(d^2,p_3,p_4,p_5,p_6)$ where $p_3+p_4+p_5+p_6=4d-6$, and since it has to be dominated by the generic Jordan type we must have $p_i\leq d-1$. The only possible partition in this case is $\bigl(d^2,(d-1)^2,(d-2)^2\bigr)$.

Let's now consider linear forms $\ell=a_0x+a_1y+a_2z+b_0u+b_1v$ which are not weak Lefschetz elements. Note that we always have $(x,y,z)^2 \subseteq \Ann F$. In the case of Lemma \ref{lemma:min_perazzo}(\ref{minimalHFPerazzo1}), the linear form $\ell$ not being a weak Lefschetz element  means $b_0=0$, according to Proposition \ref{prop:Lefschetz-elements}. As every monomial in $x,y,z,v$ of degree four belongs to $\Ann F$ we then have $\ell^4=0$ in $A_F$. Therefore it's impossible to have a Jordan string of length more than four. 

We move on to the case of Lemma \ref{lemma:min_perazzo}(\ref{minimalHFPerazzo2}). By Proposition \ref{prop:Lefschetz-elements} we have $b_0=0$ or $b_1=0$. If $b_0=b_1=0$ we have $\ell^2=0$ in $A_F$, and hence all strings have length at most two. Suppose $b_0=0$ and $b_1 \ne 0$. We may assume $b_1=1$. As $xv,yv^2 \in \Ann F$ we have
\[
\ell^s=sa_2zv^{s-1}+v^{s}, \quad \text{and} \quad \ell^{s-1}z=zv^{s-1},
\]
for any $s \ge 3$. This gives us two strings of length $d+1$ and $d-1$, and they span all polynomials in $A_F$ in only the variables $z$ and $v$. Using also the fact that $uv^2 \in \Ann F$ we get $\ell x=\ell^2y=\ell^3u=0$ in $A_F$,  and hence all remaining strings must have length at most three. Suppose instead $b_0=1$ and $b_1=0$. As $zu \in \Ann F$ we have
\[
\ell^s=sa_0xu^{s-1}+sa_1yu^{s-1}+u^s
\]
for $s \ge 2$ and
\[
\ell^{d-1}=(d-1)a_0xu^{d-2}+u^{d-1} \ne 0.
\]
Moreover 
\[
\ell^{d-2}x=xu^{d-2}, \quad \ell^{d-2}y=yu^{d-2}, \quad \ell^{d-2}v=(d-2)a_1yu^{d-3}v+u^{d-2}v.
\]
It is straightforward to verify that $\ell^{d-1}$, $\ell^{d-2}x$, $\ell^{d-2}y$, and $\ell^{d-2}v$ acting on $F$ give four linearly independent elements, and therefore this gives us four Jordan strings of length at least $d-1$. We can also see that the one-dimensional space $[A_F]_d$ is covered, as $\ell^{d-1}x=xu^{d-1} \ne 0$. 
There is one more string starting in degree one, and as $\ell z=0$ this string has length one. In degree $1<s<d$ the monomials $v^s$ and $zv^{s-1}$ span a two dimensional subspace of $[A_F]_s$. Both monomials are in the kernel of multiplication by $\ell^2$,  therefore they cannot be in the span of $\ell^s$, $\ell^{s-1}x$, $\ell^{s-1}y$, $\ell^{s-1}v$, except possibly when $s=d-2$. In any case this shows that the remaining strings have length at most two.

Finally we consider the last  case of Lemma \ref{lemma:min_perazzo}(\ref{minimalHFPerazzo3}). Here $\ell$ is non-Lefschetz if $b_0=0$, $b_1=0$, or $b_0+\lambda b_1=0$. Let's start by assuming $b_0=0$ and $b_1=1$. In this case 
\[
\ell^s=sa_1yv^{s-1}+sa_2zv^{s-1}+v^s \quad \text{when} \ s \ge 2.
\]
The power $\ell^{d-1}$ together with 
\[
\ell^{d-2}y=v^{d-2}y, \quad \ell^{d-2}z=v^{d-2}z, \quad \ell^{d-2}u=(d-2)a_1yuv^{d-3}+uv^{d-2}  
\]
are linearly independent, and give us four strings of lengths at least $d-1$. As $\ell x=0$ the remaining string starting in degree one will have length one. In degree $1<s<d$ the elements $u^{s-1}(\lambda u-v)$ and $xu^{s-1}$ are linearly independent, and both in the kernel of multiplication by $\ell^2$. Hence all remaining strings have lengths at most two. 

The case $b_0=1$ and $b_1=0$ is treated analogously. Suppose $b_0+\lambda b_1=0$. We may assume $b_0=\lambda$ and $b_1=-1$. Note that $y(\lambda u -v) \in \Ann F$ in this case. We have 
\[
\ell^s=s \lambda a_0 xu^{s-1}-s a_2 zv^{s-1} + (\lambda u -v)^s \quad \text{for} \ s \ge 2.
\]
It is straightforward to check that $\ell^{d-1}$ together with 
\begin{align*}
& \ell^{d-2}x=x(\lambda u -v)^{d-2}, \quad \ell^{d-2}z=z(\lambda u -v)^{d-2},\\
 & \text{and} \quad \ell^{d-2}u=(d-2) \lambda a_0 xu^{d-2}+(\lambda u -v)^{d-2}u
\end{align*} 
form a $4$-dimensional space. Hence the Jordan basis has four strings of length at least $d-1$. As $\ell y=0$ in $A_F$ we get one string of length one starting in degree one. In higher degree the elements $y(u+\lambda v)^{s}$ and $u^sv$ form a two-dimensional space contained in the kernel of multiplication by $\ell^2$. It follows that all remaining strings will have length at most two. 
\end{proof}

\subsection{Jordan type of a non-Lefschetz element}
We do a more thorough study of the possible Jordan types in the case of Lemma \ref{lemma:min_perazzo} (\ref{minimalHFPerazzo1}). As the possible Jordan types of Lefschetz elements were discussed in detail above, let now $\ell=a_0x+a_1y+a_2z+b_0u+b_1v$ be a non-Lefschetz element. We present the Jordan types, as well as the strings in a pre\-\mbox{-Jordan} basis, in each case.


\begin{center}
    \begin{tabular}{|c|c|l|}
        \hline
        Conditions & $P_{\ell,A_F}$ & pre-Jordan basis\\
        \hline
        \hline
        \multirow{4}{4em}{$a_0\in K$, $a_1\in K$, $a_2\ne 0$, $b_1\ne 0$}  & \multirow{4}{5,5em}{$(4^{d-2},2^d,1^2)$} & $u^i\mapsto \ell u^i\mapsto \ell^2 u^i \mapsto \ell^3 u^i,\, 0\le i\le d-3$\\
        &&$u^iv\mapsto \ell u^iv,\, 0\le i\le d-3$\\
        &&$y\mapsto \ell y$, $u^{d-2}\mapsto \ell u^{d-2}$ \\
        &&$x,u^{d-1}$\\
        \hline
        
        \multirow{4}{4em}{$a_0\in K$, $a_1\in K$, $a_2=0$, $b_1\ne 0$}  & \multirow{4}{5,5em}{$(3^{2d-4},2^{2},1^2)$} & $u^i\mapsto \ell u^i\mapsto \ell^2 u^i,\, 0\le i\le d-3$\\
        && $zu^i\mapsto \ell zu^i \mapsto \ell^2 zu^i,\, 0\le i\le d-3$\\
        && $y\mapsto \ell y$, $u^{d-2}\mapsto \ell u^{d-2}$\\
        &&$x,u^{d-1}$\\
        \hline

        \multirow{4}{4em}{$a_0\in K$, $a_1\in K$, $a_2\ne 0$, $b_1=0$}  & \multirow{4}{5,5em}{$(2^{3d-6},1^6)$} & $u^i\mapsto \ell u^i,\, 0\le i\le d-3$\\
        && $u^iv\mapsto \ell u^iv,\, 0\le i\le d-3$\\
        && $u^iv^2\mapsto \ell u^iv^2,\, 0\le i\le d-3$\\
        &&$x,y,xu,u^{d-2},u^{d-2}v,u^{d-1}$\\
        \hline


        \multirow{4}{4em}{$a_0\in K$, $a_1\ne 0$, $a_2=0$, $b_1=0$}  & \multirow{4}{5,5em}{$(2^{2d-2},1^{2d-2})$} & $u^i\mapsto \ell u^i,\, 0\le i\le d-2$\\
        &&$u^iv\mapsto \ell u^iv,\, 0\le i\le d-2$\\
        &&$zu^i,u^iv^2,\, 0\le i \le d-3$ \\
        &&$x, u^{d-1}$\\
        \hline
        
        \multirow{4}{4em}{$a_0\ne 0$, $a_1=0$, $a_2=0$, $b_1=0$}  & \multirow{4}{5,5em}{$(2^{d},1^{4d-6})$} & \multirow{4}{10em}{$u^i\mapsto \ell u^i,\, 0\le i\le d-1$ $u^iv^2,zu^i,\, 0\le i \le d-3$ $u^iv,yu^i,\, 0\le i \le d-2$}\\
        &&\\
        &&\\
        &&\\
        \hline
    \end{tabular}
\end{center}

We notice that the computations agree with Theorem \ref{thm:JT}. We also remark that these Jordan types are in a chain with respect to the dominance order:
\[
(4^{d-2},2^d,1^2)>(3^{2d-4},2^2,1^2)>
  (2^{3d-6},1^6)>(2^{2d-2},1^{2d-2})>
  (2^{d},1^{4d-6}).
\]


\newcommand{\etalchar}[1]{$^{#1}$}
\providecommand{\bysame}{\leavevmode\hbox to3em{\hrulefill}\thinspace}
\providecommand{\MR}{\relax\ifhmode\unskip\space\fi MR }
\providecommand{\MRhref}[2]{%
  \href{http://www.ams.org/mathscinet-getitem?mr=#1}{#2}
}
\providecommand{\href}[2]{#2}

\end{document}